\documentclass[preprint,12pt]{elsarticle}
\graphicspath{ {./figures/} }
\usepackage{subfig}
\usepackage{xcolor}
\usepackage{hyperref}

\usepackage{amsmath,amsthm}
\usepackage{amsfonts}
\usepackage{amssymb}
\usepackage{textcomp}
\usepackage{stmaryrd}
\usepackage{tikz}
\usepackage{breqn}
\usepackage{dsfont}
\usepackage{pdfpages}
\usepackage{graphicx}
\usepackage{tipa}
\usepackage{upgreek}

\newtheorem{Definition}{Definition}

\newtheorem{Theorem}{Theorem}

\journal{}
\begin{document}
\begin{frontmatter}
\title{Can a Fractional Order Delay Differential Equation be Chaotic Whose Integer-Order Counterpart is Stable? }
\author{Sachin Bhalekar\footnote{Corresponding Author Email: sachinbhalekar@uohyd.ac.in}, Deepa Gupta}
 \address{School of Mathematics and Statistics, University of Hyderabad, Hyderabad, 500046 India}
\begin{abstract}
For the fractional order systems \[D^\alpha x(t)=f(x),\quad 0<\alpha\leq 1,\]  one can have a critical value of $\alpha$ viz $\alpha_*$ such that the system is stable for $0<\alpha<\alpha_*$ and unstable for $\alpha_*<\alpha\leq 1$. In general, if such system is stable for some $\alpha_0\in(0,1)$ then it remains stable for all $\alpha<\alpha_0.$\\
In this paper, we show that there are some delay differential equations \[D^\alpha x(t)=f(x(t),x(t-\tau))\] of the fractional order which behave in an exactly opposite way. These systems are unstable for higher values of fractional order and stable for the lower values.\\
The striking observation is the example which is chaotic for $\alpha=0.27$ but stable for $\alpha=1$. This cannot be observed in the fractional differential equations (FDEs) without delay.\\
We provide the complete bifurcation scenarios in the scalar FDEs.
\end{abstract}
\end{frontmatter} 
\section{Introduction}
Various equations arising in the theory of Applied Mathematics played a vital role in modeling natural systems \cite{rafikov2008mathematical,haberman1998mathematical,logan2013applied,morse1954methods,manninen2006developing,braun1983differential}. The time variable is continuous in most of these systems and one can use the ordinary or partial differential equations to model the dynamics. The rate of change of the physical quantity under consideration is described by the ``derivatives'' in these models. The remarkable results in the history of Science viz. Newtons's laws \cite{o1990newton,newton1833philosophiae}, Einstein's theory \cite{eriksson1996computational,einstein1905photoelectric}, heat conduction \cite{douglas1956numerical}, chemical kinetics \cite{steinfeld1999chemical} etc made use of the derivative operators.\\
 In spite of a wide range of applicability, the derivative operator fails in modeling the ``memory and hereditary properties'' in the natural systems because it is a ``local operator''. Vito Volterra in 1959 introduced the theory of non-local operators \cite{volterra1959theory} which proved advantageous in modeling memory and hereditary properties.\\
 The fractional derivative(FD) \cite{podlubny1998fractional,kilbas2006theory} and the delay \cite{lakshmanan2011dynamics,smith2011introduction} are the popularly used mathematical concepts which are very useful in achieving this goal. The fractional differential equation can be viewed as an integral equation with the kernel $(t-\tau)^{\alpha-1}/\Gamma(\alpha)$ (or similar). The flexible order $\alpha$ of the derivative can be adjusted according to the physical system to achieve more accuracy in the modeling. On the other hand, the delay differential equation(DDE) contains the value of state variable at the past time to include the memory of the system \cite{bocharov2000numerical,gopalsamy2013stability}.\\
 Fractional differential equations have a wide applications in the studies of viscoelastic materials \cite{mainardi2010fractional,torvik1984appearance} as well as in many fields of science and engineering including fluid flow, rheology, electrical networks and electromagnetic theory \cite{debnath2003recent,kilbas2006theory,saqib2019application}. Even a non-linear oscillation of an earthquake can be modeled with the fractional derivatives \cite{koh1990application}. These operators are also proved useful in modeling various phenomena in bioenginering and related area \cite{magin2004fractional,magin2008modeling}.\\
 The delay appears in the natural system because there is always a time duration of certain hidden processes (e.g. duration of infectious period) \cite{rihan2014time}. The delay differential equations have applications in various areas of life sciences such as in population dynamics, epidemiology, immunology, physiology and neural networks \cite{bocharov2000numerical,rihan2014delay,baker1998modelling}. So, the fractional order delay differential equations (FDDE) containing FD as well as the delay seem to be more realistic.\\ 
 Although the fractional derivative and delay are more practical, the problems become infinite dimensional \cite{smith2011introduction,garrappa2020initial}.
 The mathematical analysis is more involved and one has to be careful while applying the classical results.\\
 The stability and bifurcation analyses are the principal branches of Dynamical Systems \cite{meiss2007differential,hale2012dynamics}. One has to provide these results before applying these useful tools (viz. FD and delay) in the real-life applications. In the seminal work, Bhalekar \cite{bhalekar2016stability} provided the stability and bifurcation analyses of the FDDE system\[D^\alpha x(t)=ax(t)+bx(t-\tau).\] 
 He provided three regions viz. unstable(U), stable(S) and the delay dependent stable(DS) regions (cf. Figure \ref{fig1.2}).\\
 In the delay dependent stable region, the system is stable for some values of delay $\tau$ and become unstable for other values.\\
 In this work, we extend the analysis in \cite{bhalekar2016stability} and provide more details on the bifurcation regions in the delay dependent stable region.\\
 The manuscript is organized as below: Section \ref{sec3} provides the basic definitions and results useful in understanding this paper. Section \ref{sec2} analyses the stability and bifurcations. Section \ref{sec4} contains illustrative examples. In the Section \ref{sec5}, we provide some new observations. The conclusions are given in Section \ref{sec6}.
\section{Preliminaries}\label{sec3}
In this section, we provide some basic definitions described in the literature \cite{podlubny1998fractional,lakshmanan2011dynamics,kilbas2006theory,bhalekar2013stability,diethelm2002analysis}.\\
\begin{Definition}[Fractional Integral]
For any $f \in \mathcal{L}^{1}(0,b)$ the Riemann-Liouville fractional integral of order $\upmu >0$, is given by 

\begin{equation*}
\textit{I}^\upmu f(t)=\dfrac{1}{\Gamma(\upmu)}\int_{0}^{t}(t-\tau)^{\upmu-1}f(\tau)d\tau  , \quad   0<t<b.
\end{equation*}

\end{Definition}
\begin{Definition}[Caputo Fractional Derivative]
For $f \in \mathcal{L}^{1}(0,b)$, $0<t<b$ and $m-1<\upmu\leq m$, $m \in \mathbb{N}$, the Caputo fractional derivative of function $f$ of order $\upmu$ is defined by,
\[\textit{D}^{\upmu} f(t)=
\begin{cases}
\frac{d^m}{dt^m} f(t) ,\textit{ if } \quad \upmu = m \\ \textit{I}^{m-\upmu}\dfrac{d^m f(t)}{dt^m}, \textit{ if } \quad  m-1< \upmu < m. 
\end{cases}\]
 Note that for  $m-1 < \upmu \leq m$, $m\in \mathbb{N},$

\[\textit{I}^\upmu\textit{D}^\upmu f(t)=f(t)-\sum_{k=0}^{m-1}\dfrac{d^k f(0)}{dt^k}\dfrac{t^k}{k!}.\]

\end{Definition}
\begin{Definition}[Equilibrium Point]
 Consider the generalized delay differential equation
  \begin{equation}
  D^\alpha{x}(t)=f(x(t),x(t-\tau)),  \quad 0<\alpha\leq 1,\label{eq1}
  \end{equation}
  where $\tau \textgreater0$,  $f: E\rightarrow \mathbb{R}$, $E\subseteq \mathbb{R}^2$ is open and $f\in C^1(E)$.

  A steady state solution of equation (\ref{eq1}) is called an equilibrium point.
  
  Note that $x^*$ is an equilibrium point if and only if
  \begin{equation}
  f(x^*,x^*)=0.
  \end{equation}
  \end{Definition}
  The linearization of equation (\ref{eq1}) near $x^*$ is given by \cite{bhalekar2016stability}.
  \begin{equation}
   D^{\alpha} x(t)=ax(t)+b x(t-\tau),\label{eq3}
   \end{equation} 
   where $a=\partial_1f(x^*,x^*)$ and $b=\partial_2f(x^*,x^*)$ are partial derivatives of $f$ with respect to the first and second variables respectively evaluated at $x(t)=x(t-\tau)=x^*.$ 
  Near the equilibrium point, the qualitative behaviour of the trajectories of equations (\ref{eq3}) and (\ref{eq1}) are same.
  \section{Bifurcation Analysis}\label{sec2}
\begin{figure}
 \centering
\includegraphics[width=0.6\textwidth]{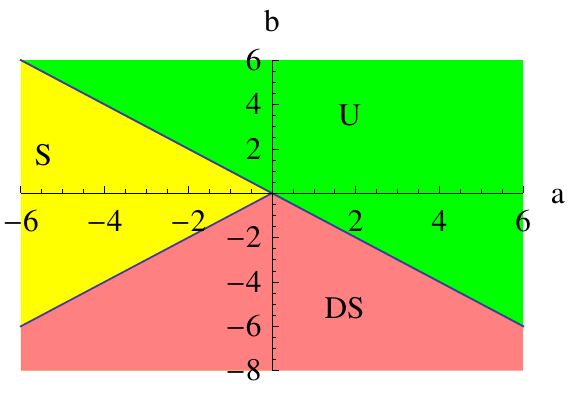}
\caption{ Stability regions of equation (\ref{eq3})  }\label{fig1.2}
\end{figure}

\begin{figure}
 \centering
\includegraphics[width=0.9\textwidth]{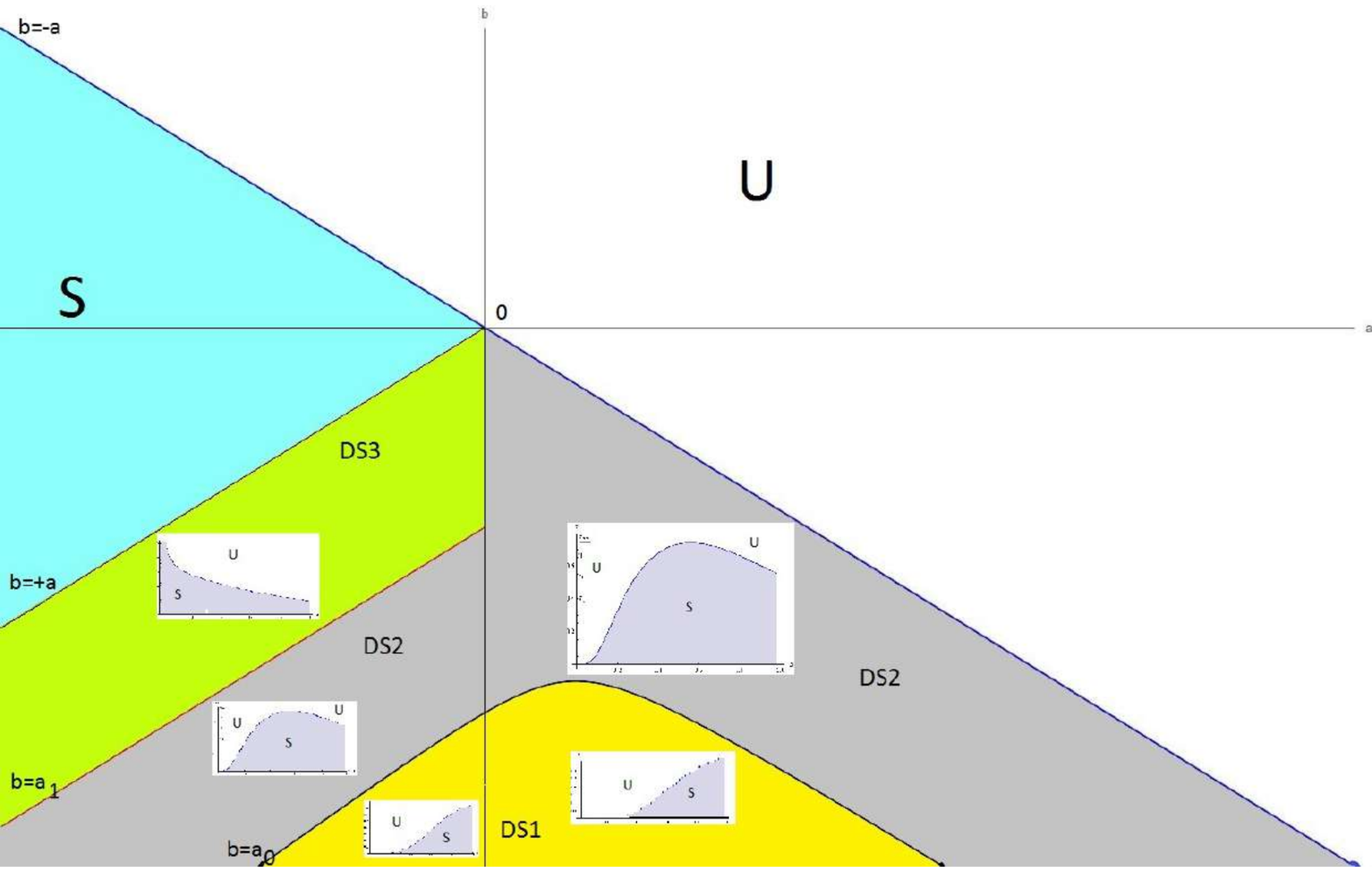}
\caption{The complete bifurcation diagram of equation (\ref{eq3})}\label{fig1.1}
\end{figure}

  \begin{Theorem}\cite{bhalekar2016stability}\label{thm1}
   Consider the fractional order linear delay differential equation (\ref{eq3}).\\
\textbf{Case 1} If b $\in (-\infty,-|a|)$ then the stability region of $x^*$ in $(\tau,a,b)$ parameter space is located between the plane $\tau=0$ and
\begin{equation}\label{eq1.8}
 \tau(a,b,\alpha)=\dfrac{\arccos\Bigg(\dfrac{\Bigg(acos\Big(\dfrac{\alpha\pi}{2}\Big)+\sqrt{b^2-a^2\sin^2\Big(\dfrac{\alpha\pi}{2}}\Big)\Bigg)\cos\dfrac{\alpha\pi}{2}-a}{b}\Bigg)}{\Bigg(a\cos\Big(\dfrac{\alpha\pi}{2}\Big)+\sqrt{b^2-a^2\sin^2\Big(\dfrac{\alpha\pi}{2}\Big)}\Bigg)^{1/\alpha}}.
 \end{equation}
The equation undergoes Hopf bifurcation at this value.\\
 \textbf{Case 2} If $b\in (-a,\infty)$ then $x^*$ is unstable for any $\tau\geq 0.$\\
 \textbf{Case 3} If $b\in (a,-a)$ and $a<0$ then $x^*$ is stable for any $\tau\geq 0.$\\
 The curve $\tau(a,b,\alpha)$ given by equation (\ref{eq1.8}) is called the boundary curve.\\ 
 \end{Theorem}
 This theorem provides three regions in $ab-$plane (cf. Figure \ref{fig1.2} ).\\
 viz. S- stable region-Case 3.\\
      U- unstable region-Case 2 and\\
      DS- delay dependent stable region-Case 1.\\
 In this work, we provide more details of the region DS in Figure \ref{fig1.2}.\\ 
 The bifurcation regions are shown in Figure \ref{fig1.1} and can be described as below:\\
\textbf{1.} U: If $a+b>0$ then system (\ref{eq3}) is unstable $\forall \tau\geq 0.$\\
\textbf{2.} S: If $a<0$ and $a<b<-a$ then the system (\ref{eq3}) is asymptotically stable $\forall \tau\geq 0.$\\
\textbf{3.} DS1: In this region, we get the stable region in $\alpha\tau-$ plane as shown in Figure \ref{fig15}.\\
\textbf{4.} DS2: The stability properties in this region are as shown in Figure \ref{fig17}.\\
\textbf{5.} DS3: The stable region in $\alpha\tau$-plane for this case is shown in Figure \ref{fig18}.\\

 \subsection{Bifurcation Scenarios}
\begin{figure}
 \centering
\includegraphics[width=0.6\textwidth]{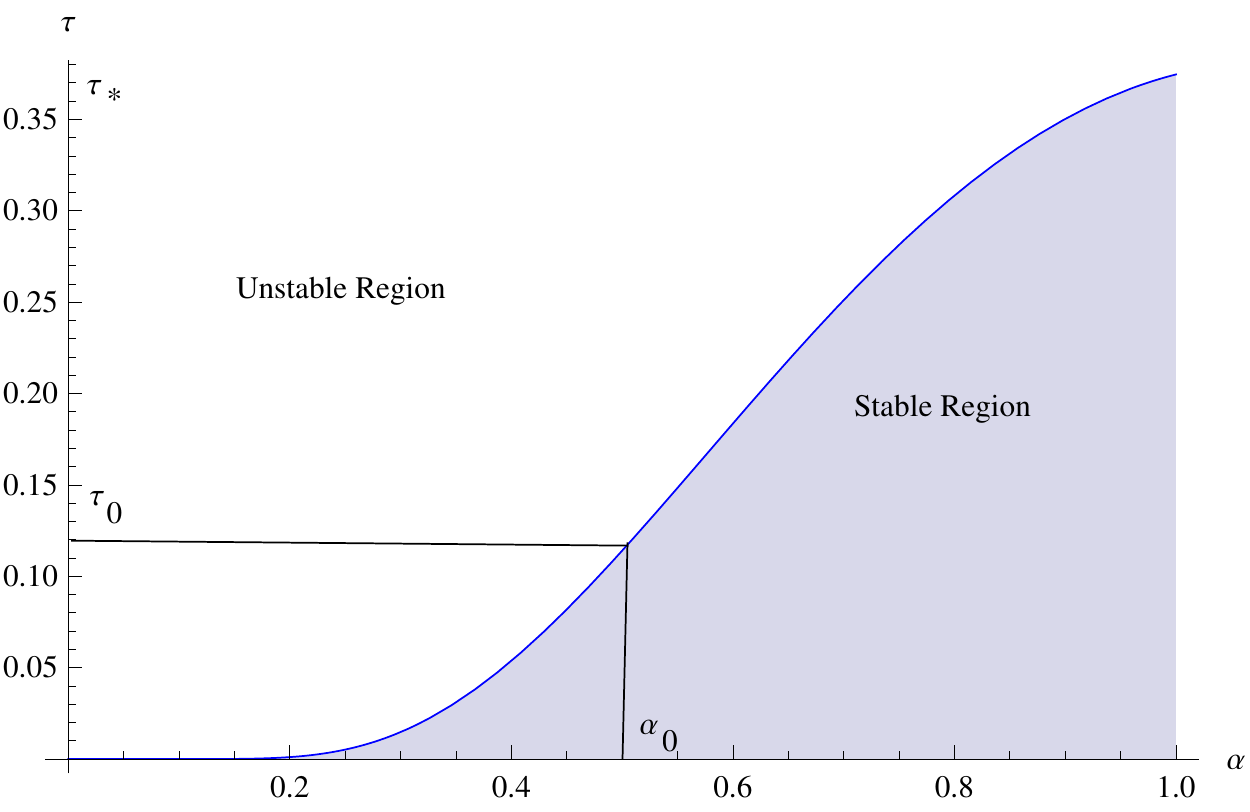}
\caption{ Boundary curve (\ref{eq1.8}) in the region DS1  }\label{fig15}
\end{figure}

In this section, we provide the details of various stable regions described in previous section.\\
\textbf{Region DS1}: It is observed that, for given $a\in\mathds{R}$, $\exists$ a function $a_0(a)$ such that $-\infty<b<a_0(a)\Rightarrow$ the boundary curve $\tau(a,b,\alpha)$ in the $\alpha\tau-$ plane is monotonically increasing as given in Figure \ref{fig15}. Let $\tau_*=\tau(a,b,1)$. For fixed $\tau_0\in(0,\tau_*]$, $\exists$ critical value of $\alpha$, say $\alpha_0$ such that equilibrium point is asymptotically stable for $\alpha_0<\alpha\leqslant1$ and unstable for $0<\alpha<\alpha_0$. If $\tau>\tau_*$ then the system is unstable for all $\alpha\in(0,1]$.\\
\textbf{Note} This case is exactly opposite of non-delayed case.\\
\subsubsection{Estimation for the function $a_0(a)$:}
We observed that for some values of $a\in\mathds{R}$, there exists functions $a_0(a)$ and $a_1(a)$ such that the curve $\tau(a,b,\alpha)$ is monotonically increasing for $b\in(-\infty,a_0(a))$ and in the interval $a_0(a)<b<a_1(a)$, the curve has local maxima. This local maxima occurs at $\alpha=1$ if $b \shortrightarrow a_0$.\\
Therefore, $a_0$ is the function such that-\\
\[\lim_{b\shortrightarrow a_0}\frac{\partial\tau(a,b,\alpha)}{\partial \alpha}\mid_{\alpha=1}=0.\]
Solving this equation, we get,
\begin{equation*}
-\arccos(\frac{-a}{a_0})=\frac{\sqrt{a_0^2-a^2}(\frac{\pi}{2})}{a\frac{\pi}{2}+\sqrt{a_0^2-a^2}\log(\sqrt{a_0^2-a^2})}
\end{equation*}
or,
\begin{equation}
\frac{-a}{a_0}=\cos\Big[\frac{-\sqrt{a_0^2-a^2}(\frac{\pi}{2})}{a\frac{\pi}{2}+\sqrt{a_0^2-a^2}log(\sqrt{a_0^2-a^2})}\Big].\label{eq2}
\end{equation}
Equation (\ref{eq2}) is the implicit expression for the function $a_0(a)$.\\
\textbf{ Region DS2: }\\
\textbf{Case 1} We assume that $a<0$.\\ If $a_0<b<a_1$ then the stable curve $\tau(a,b,\alpha)$ is not monotonic (cf. Figure(\ref{fig17})) and $\exists$ a local maxima at some $\alpha_{**}<1$. In this case there exists $\tau_*=\tau(a,b,1)$ and $\tau_{**}=\tau(a,b,\alpha_{**})$ such that-\\
$\bullet$ If $0<\tau<\tau_*$ then-\\
(i) there exists critical value of $\alpha$, say $\alpha_0$ such that $0<\alpha<\alpha_0$, we get unstable solution and for $\alpha_0<\alpha\leq 1$ we get stable solution.\\
$\bullet$ If $\tau_*<\tau<\tau_{**}$ then-\\
(i) there exist $\alpha_1$ and $\alpha_2$ such that for $0<\alpha<\alpha_1$ and $\alpha_2<\alpha\leq 1$ we get unstable solution and for $\alpha_1<\alpha<\alpha_2$ we get stable solution (stable window).\\
$\bullet$ If $\tau>\tau_{**}$ then equilibrium is unstable for all $\alpha\in(0,1].$\\
\subsubsection{Estimation of function $a_1(a):$}
\begin{figure}
 \centering
\includegraphics[width=0.6\textwidth]{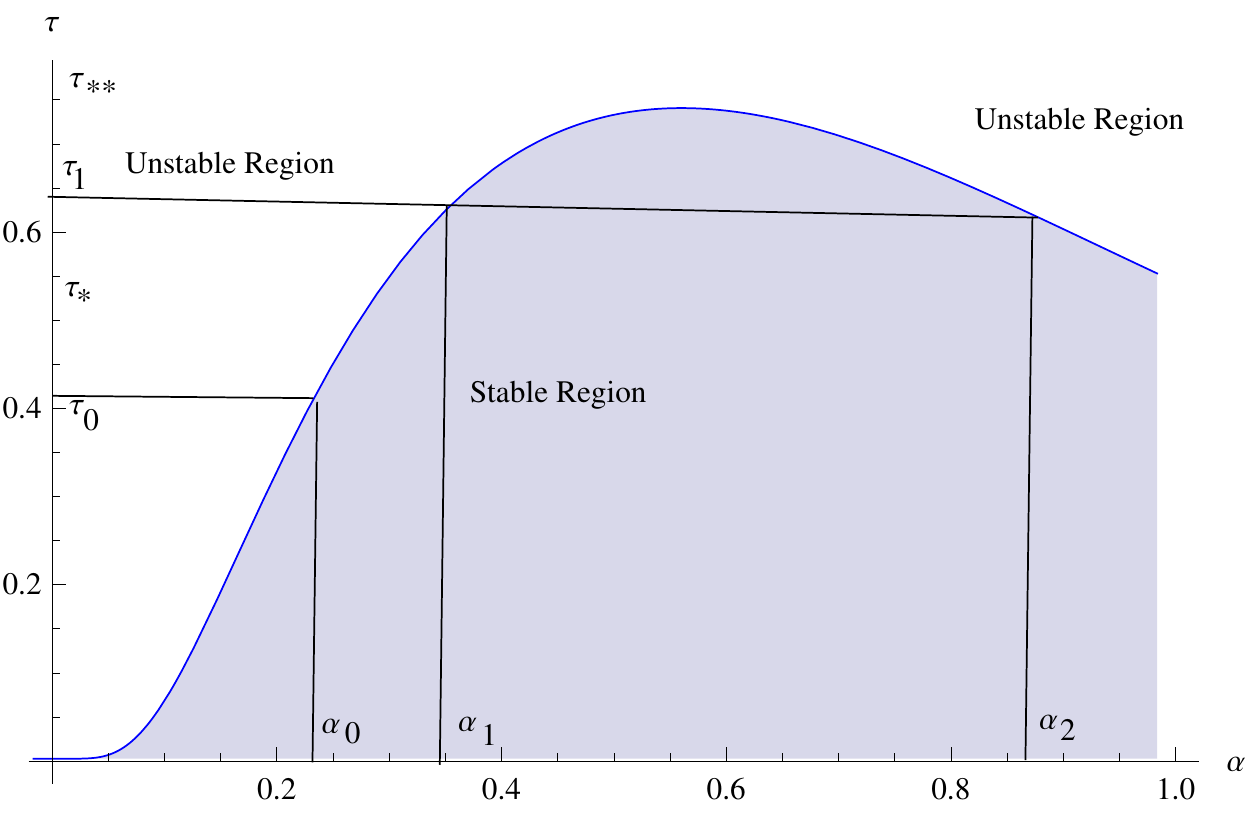}
\caption{ Boundary curve (\ref{eq1.8}) in the region DS2  }\label{fig17}
\end{figure}

\begin{figure}
 \centering
\includegraphics[width=0.6\textwidth]{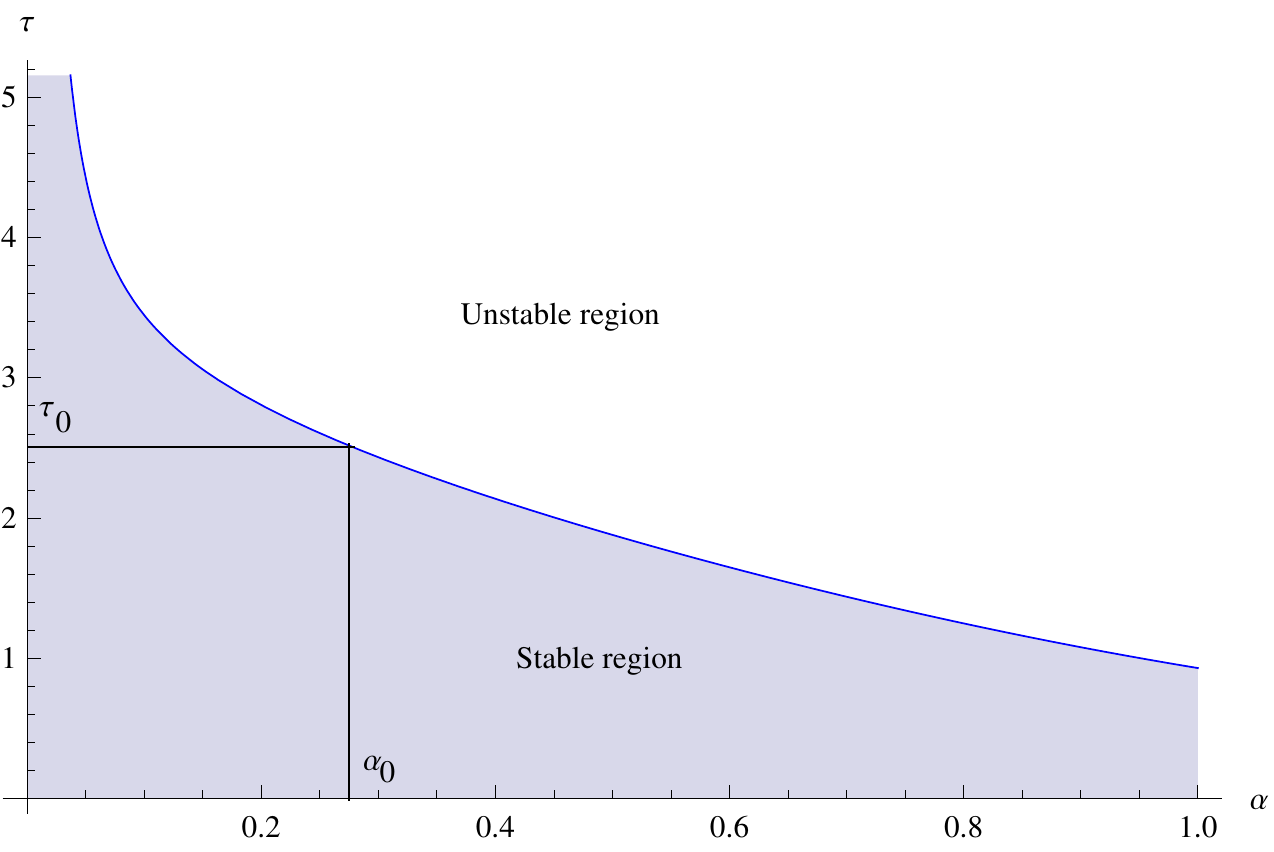}
\caption{ Boundary curve (\ref{eq1.8}) in the region DS3 }\label{fig18}
\end{figure}

Let $a<0$. As we take the values of $b>a_0$, the left end of critical curve i.e. $\tau(a,b,0)$ is zero upto some value $a_1(a)$. For $b>a_1(a)$, $\lim_{\alpha\shortrightarrow 0}\tau(a,b,\alpha)=\infty$. It is observed that, at this particular value $b=a_1(a),$\\
\[\tau(a,a_1(a),0)=\pi.\]
Solving this equation, we get ${a_1(a)=a-1}$.\\
\textbf{Case 2} We assume that $a>0.$\\ It is observed that the stable region in $\alpha\tau$-plane is shown in Figure \ref{fig17} for $a_0(a)<b<-a.$\\
Note that, $\nexists$ $a_1(a)$ if $a>0.$\\
\textbf{Region DS3:} This region exists only for $a<0$.\\
 If $a_1(a)<b<a<0$ then stable curve $\tau(a,b,\alpha)$ is monotonically decreasing in $\alpha\tau-$plane as shown in Figure \ref{fig18}. \\Note that $\lim_{\alpha\rightarrow 0}\tau(a,b,\alpha)=\infty,$ in this case.\\ For any $\tau_0>0$, $\exists$ $\alpha_0\in(0,1)$ satisfying $\tau(a,b,\alpha_0)=\tau_0$ such that $0<\alpha<\alpha_0\Rightarrow$ stable and $\alpha_0<\alpha\leq 1\Rightarrow$ unstable solutions.\\
This case is similar to the non-delayed FDEs.\\

\section{Illustrative Examples}\label{sec4}
\begin{figure}
 \centering
\includegraphics[width=0.6\textwidth]{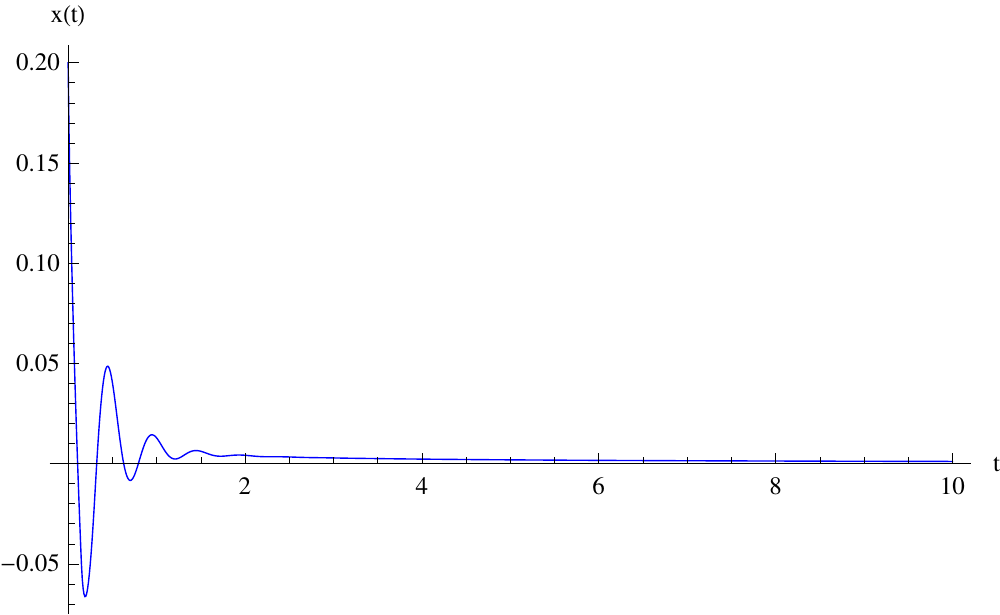}
\caption{ Stable solution of (\ref{ex1}) for $\alpha=0.8$ and $\tau=0.15$ }\label{fig20} 
\end{figure}

\begin{figure}
 \centering
\includegraphics[width=0.6\textwidth]{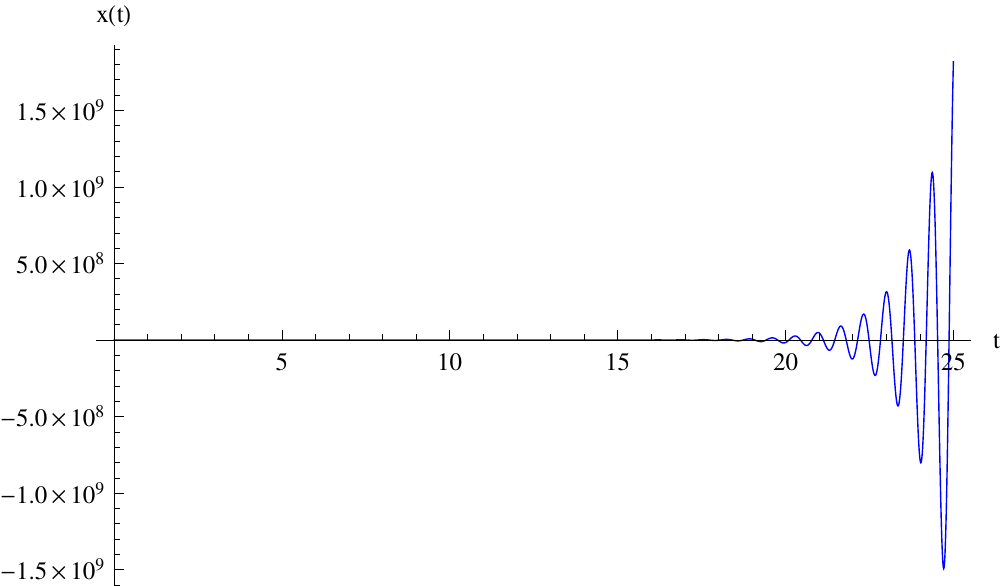}
\caption{ Unbounded solution of (\ref{ex1}) for $\alpha=0.4$ and $\tau=0.3$ }\label{fig21}
\end{figure}

\begin{figure}
 \centering
\includegraphics[width=0.8\textwidth]{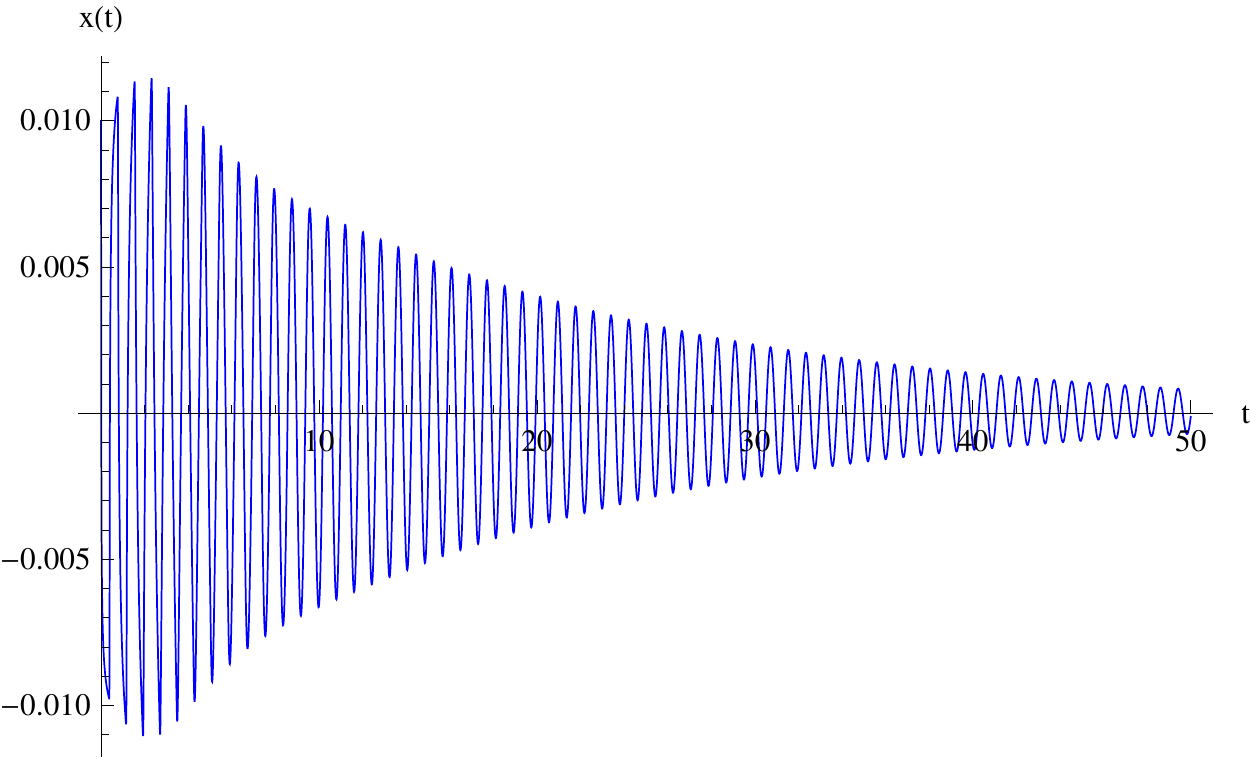}
\caption{Stable solution of equation (\ref{ex2}) for $\alpha=0.45$}\label{fig4.1}
\end{figure}

\begin{figure}
 \centering
\includegraphics[width=0.8\textwidth]{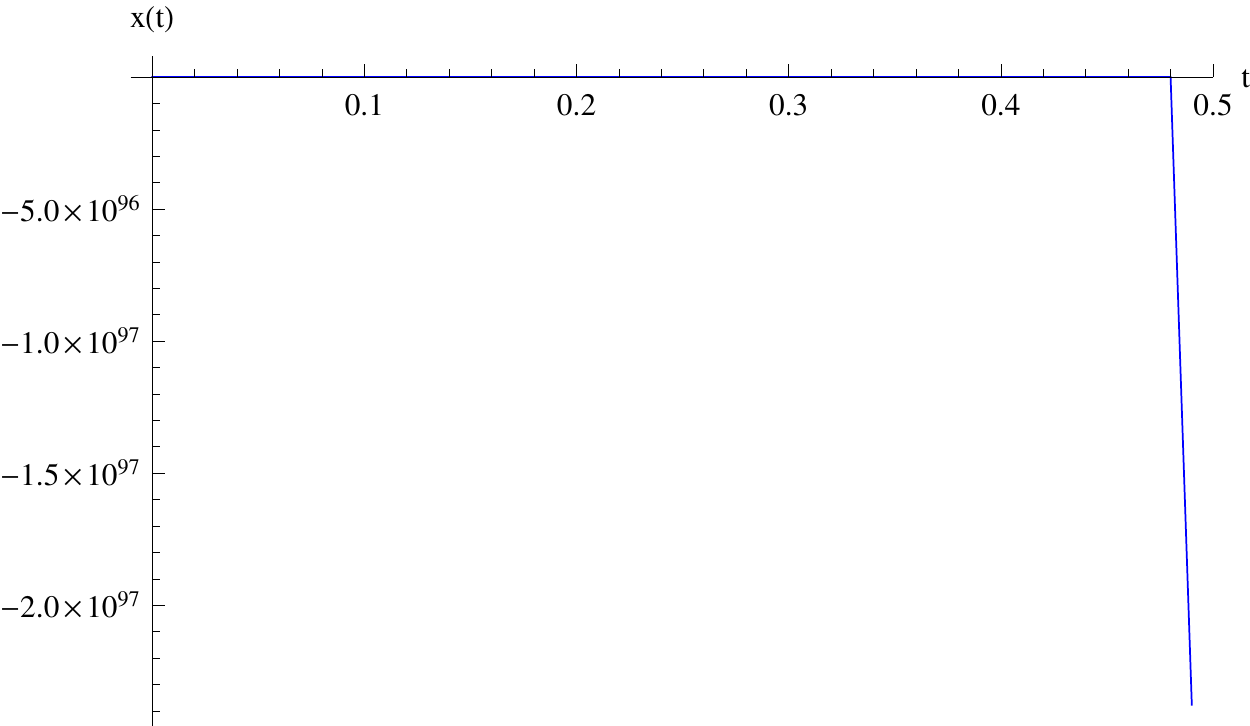}
\caption{The divergent solution of equation (\ref{ex2}) for $\alpha=0.3$ and $\tau=0.38$}\label{fig4.2}
\end{figure}

\textbf{Example 1} Consider the fractional order linear DDE
\begin{equation}
	D^\alpha x(t)=a x(t)+b x(t-\tau). \label{ex1}
\end{equation}
In this case, a=-1.5. Therefore, $a_0=-4.31582$ and $b\in(-\infty,a_0).$\\ This parameter set is in the region DS1 of Figure \ref{fig1.1}. Therefore, the stable region for this case is given by Figure \ref{fig15}.
\\The stable solution for $(\alpha,\tau)=(0.8,0.15)$ is shown in Figure \ref{fig20}. If we take $(\alpha,\tau)=(0.4,0.3)$ then it belongs to the unstable region and this is verified in Figure \ref{fig21}.\\
In the similar way we have also verified this by taking $a$ to be positive value. In this case we take $f(x(t),x(t-\tau))=2x(t)-x(t)^2-x(t-\tau)^3-4x(t-\tau)$.\\
At the equilibrium point zero we get $a=2$, $b=-4$ and $a_0=-2.33343$. By fixing $\tau=0.12$ we get critical value $\alpha=0.583977$. So, we get stable solution when $\alpha=0.7$ and for $\alpha=0.4$ we get unstable solution.\\
\textbf{ Example 2} Let us consider the non-linear system
\begin{equation}
	D^\alpha x(t)=-10.9x(t-\tau)-x(t-\tau)^3-x(t)^2-9x(t), \label{ex2}
\end{equation}
 where $0<\alpha\leq 1.$\\
Here, $f(x,y)=-9x-x^2-10.9y-y^3$. At equilibrium $x^*=0$, we have, $a=\partial_1f(0,0)=-9$ and $b=\partial_2f(0,0)=-10.9$. For this $a$, we have $a_0=-11.5972$ and $a_1=-10$. Therefore, $b\in(a_0,a_1)$.\\
Therefore, the local behaviour of $x^*$ is according to the bifurcation region DS2. The stable region will be as shown in Figure \ref{fig17} with $\tau_*=0.413437$, $\tau_{**}=0.524926$ and $\alpha_{**}=0.653835.$ \\ 
\textbf{(1)} If we fix $\tau=0.38<\tau_*$, then we get the critical value $\alpha_0=0.384137$. So, we get stable solution for $\alpha=0.45$ (cf. Figure \ref{fig4.1} ) and unstable solution for $\alpha=0.3$ (cf. Figure \ref{fig4.2}).\\

\begin{figure}
 \centering
\includegraphics[width=0.8\textwidth]{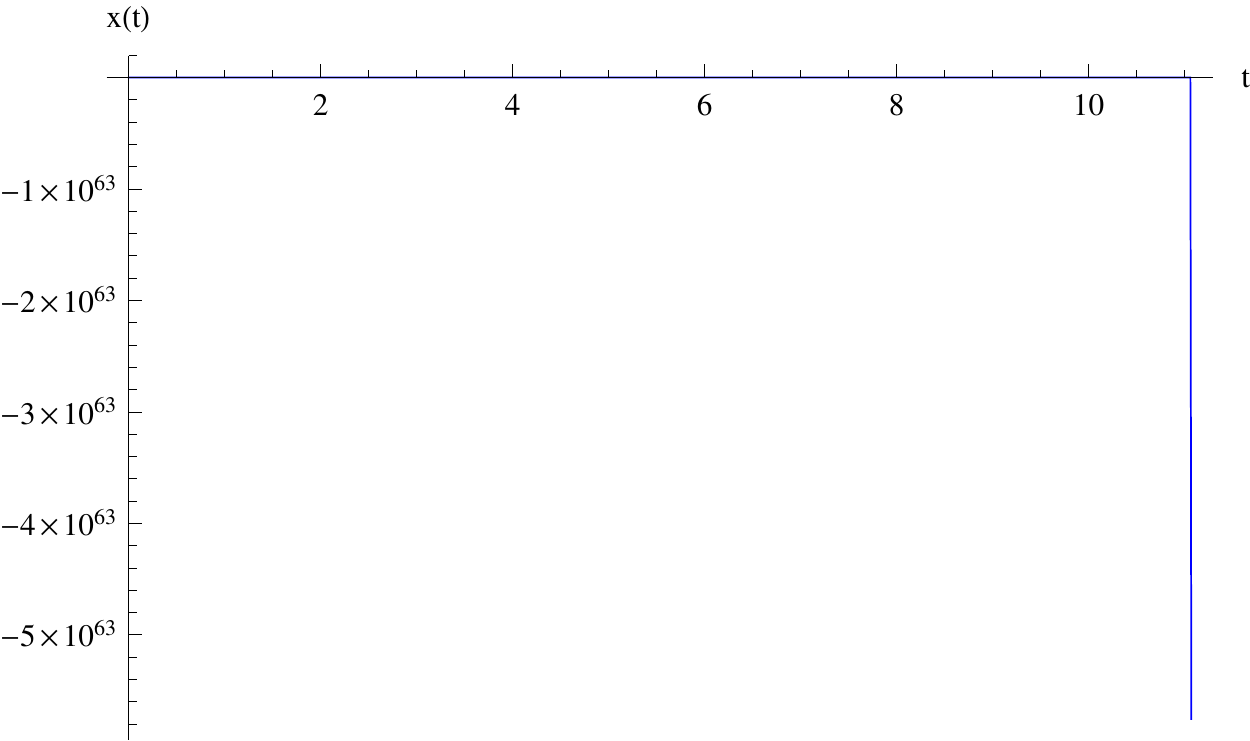}
\caption{Unstable solution of (\ref{ex2}) for $\alpha=0.41<\alpha_1$ and $\tau=0.45$}\label{fig4.3}
\end{figure}

\begin{figure}
\centering
\includegraphics[width=0.8\textwidth]{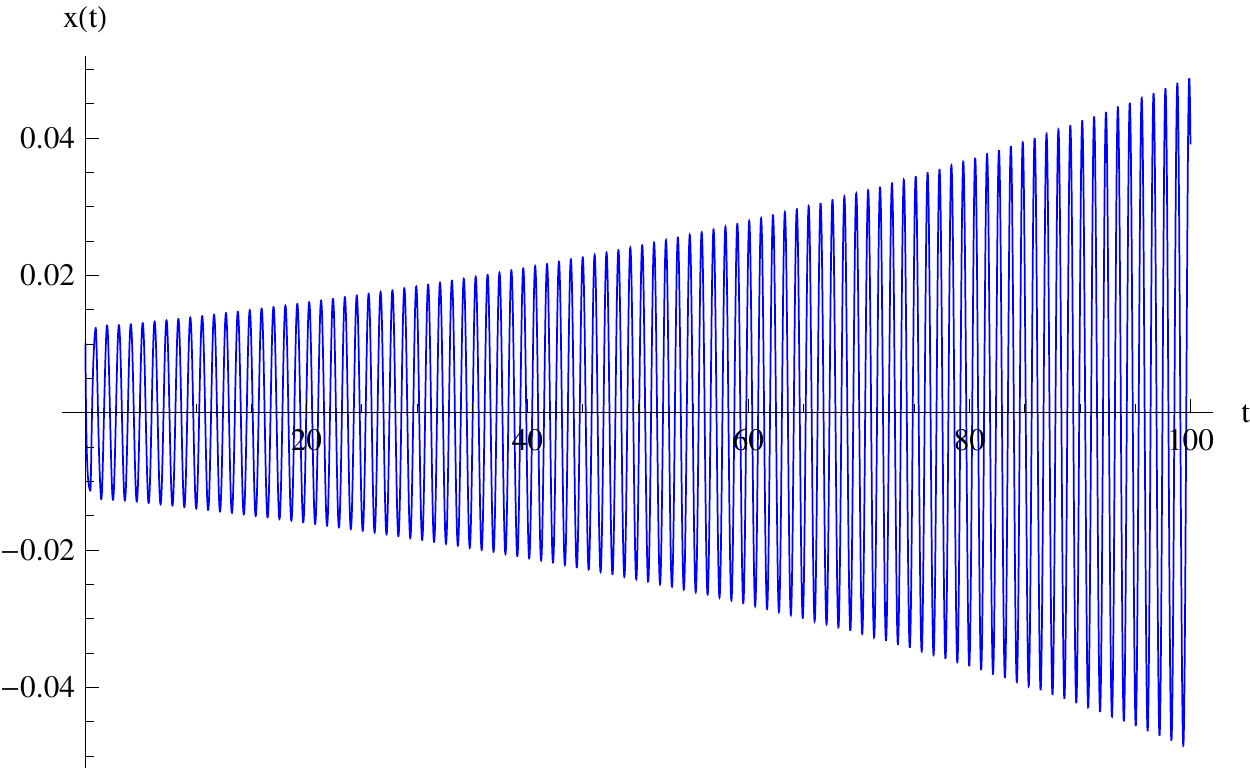}
\caption{Unstable solution of (\ref{ex2}) for $\alpha=0.95>\alpha_2$ and $\tau=0.45$}\label{fig4.4}
\end{figure}

\begin{figure}
\centering
\includegraphics[width=0.8\textwidth]{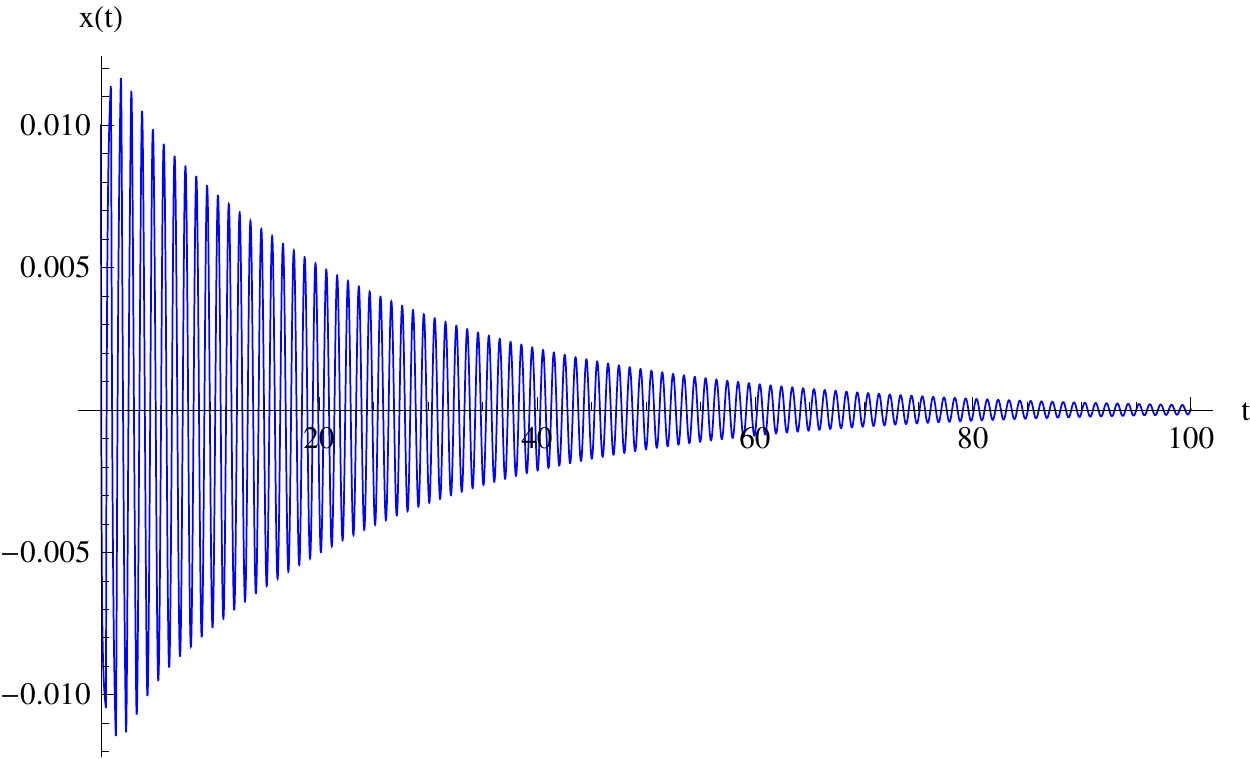}
\caption{Stable solution of (\ref{ex2}) for $\alpha=0.7$ which is in between $\alpha_1$ and $\alpha_2$ for $\tau=0.45$}\label{fig4.5}
\end{figure}

\begin{figure}
\centering
\includegraphics[width=0.8\textwidth]{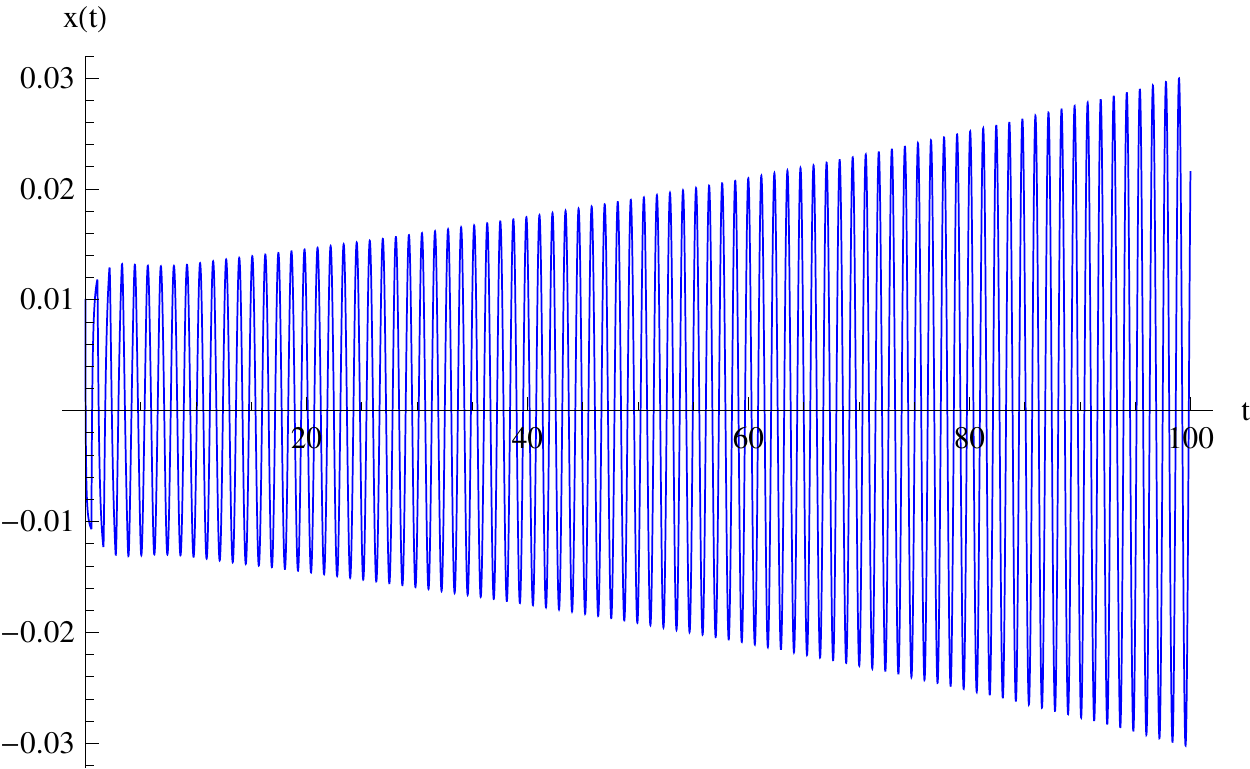}
\caption{Unstable solution of (\ref{ex2}) for $\tau=0.54>\tau_{**}$ and $\alpha=0.7$}\label{fig8.5}
\end{figure}

\textbf{(2)} When $\tau=0.45$ which is in between $\tau_*$ and $\tau_{**}$, we get $\alpha_1=0.454118$ and $\alpha_2=0.919559$. Unbounded solutions are shown in Figures \ref{fig4.3} and \ref{fig4.4} when $\alpha=0.41$ and $\alpha=0.95$ respectively and bounded solution shown in Figure \ref{fig4.5} for $\alpha=0.7.$\\
\textbf{(3)} When $\tau=0.54>\tau_{**}$ and $\alpha=0.7$ we get unstable solution (cf. Figure \ref{fig8.5}).\\ 
We also verify the results in the region DS2 with $a>0$. In this case $f(x,y)=9x-9.03y$. So, near the equilibrium point $0$, $a=9$, $b=-9.03$, $a_0=-9.06263$, $\tau_*=0.110865$ and $\tau_{**}=0.111041$. Therefore, when $\tau=0.1$ we get critical value $\alpha=0.906205$. So, we get stable solution for $\alpha=0.94$ and unstable for $\alpha=0.8$. Moreover, when $\tau=0.111$ the corresponding $\alpha_1=0.986534$ and $\alpha_2=0.99529$. So, we get divergent solutions for $\alpha=0.8$, $1$ and convergent solution for $\alpha=0.9$. Also, when we go beyond $\tau_{**}$ say $0.12$ we get unbounded solution for $\alpha=0.9.$\\

\begin{figure}
\includegraphics[width=0.8\textwidth]{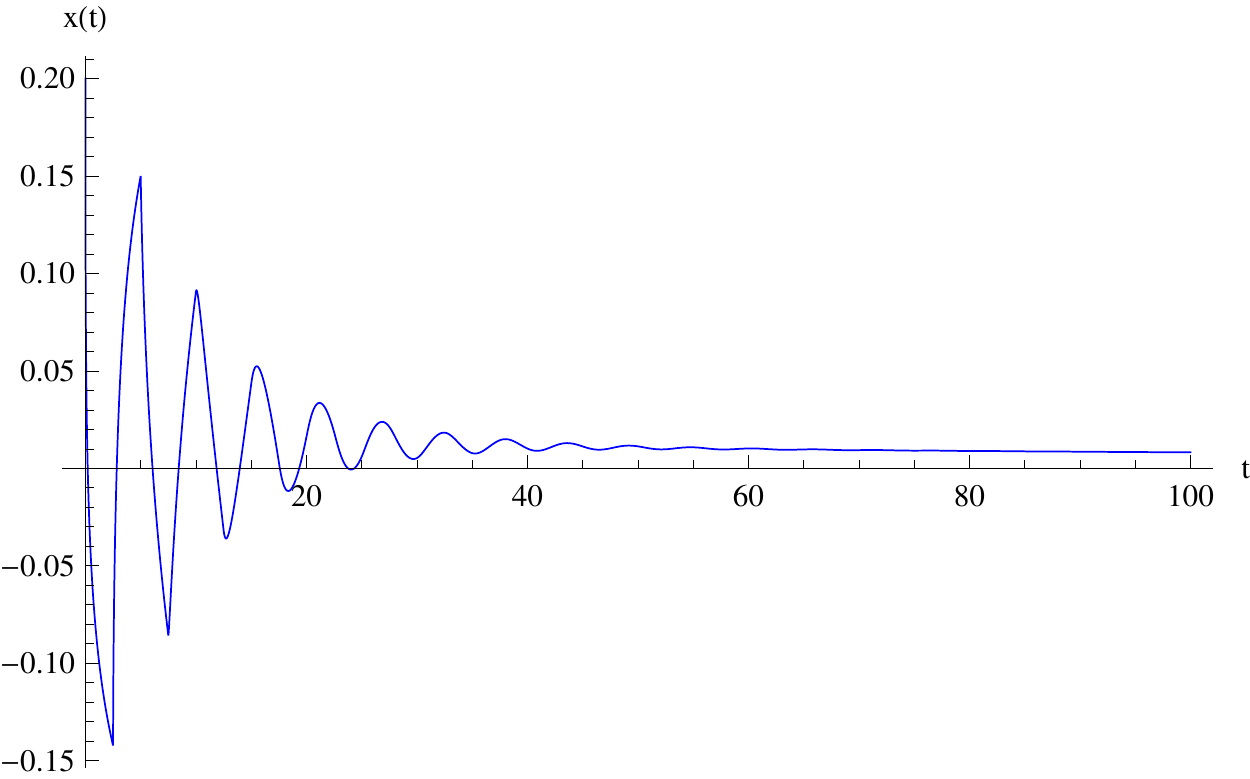}
\caption{For the Example 3 when $\alpha=0.4$ and $\tau=2.5$ we get convergent solution near equilibrium $0$}\label{fig8}
\end{figure}

\begin{figure}
\includegraphics[width=0.8\textwidth]{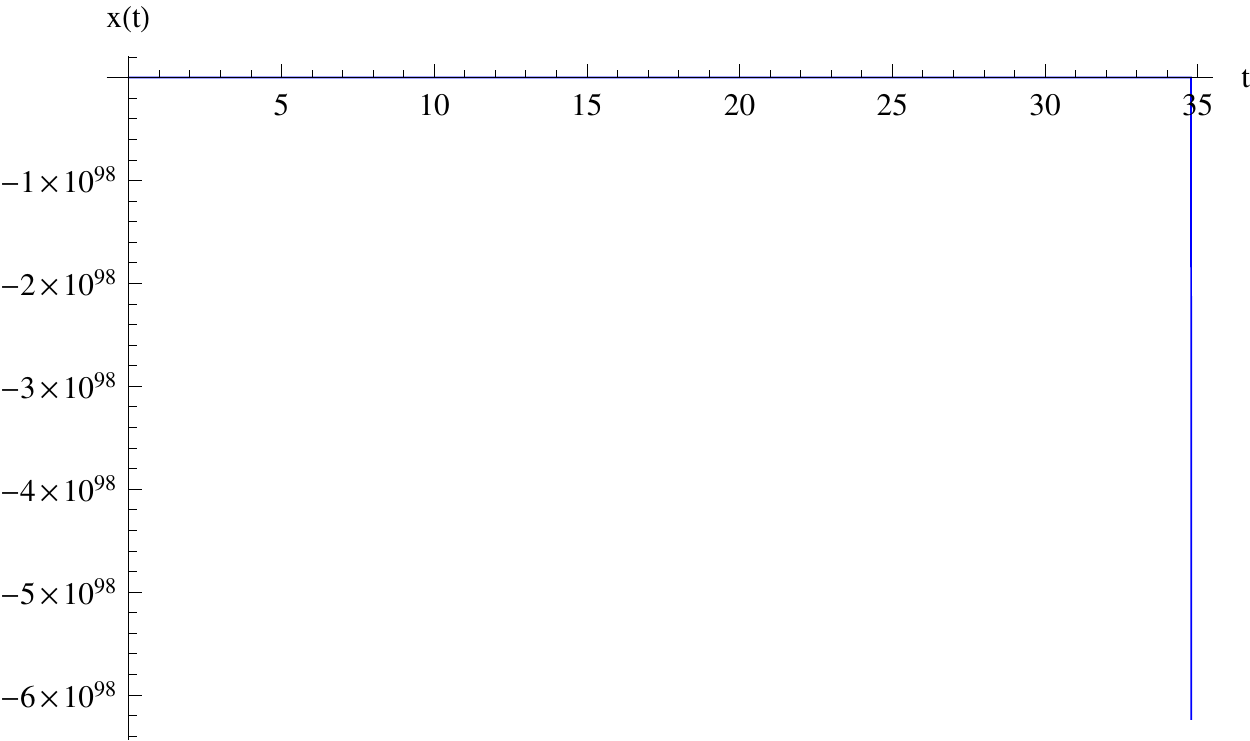}
\caption{ When $\alpha=0.91$ and $\tau=2.5$ we get divergent solution}\label{fig9}
\end{figure}

\textbf{Example 3} For the region DS3, we take a non linear FDDE (\ref{eq1}) with\\
 $f(x(t),x(t-\tau))=-x(t)-x(t)^2-1.5x(t-\tau)-x(t-\tau)^3.$\\ 
 So, near the equilibrium point $0$ we have $a=-1$ and $b=-1.5$ and $a_1=-2$. Therefore, $b\in(a_1,a)$. For the fixed $\tau_0=2.5$, we get the critical value $\alpha_0=0.904463$. Therefore, for the smaller value  $\alpha=0.4,$ we obtain stable curve (cf. Figure \ref{fig8}) and for the higher values of $\alpha$, e.g. $0.91$ we get divergent solution (cf. Figure (\ref{fig9})).

\section{New Observations}\label{sec5}
\begin{figure}
\centering
\includegraphics[width=0.6\textwidth]{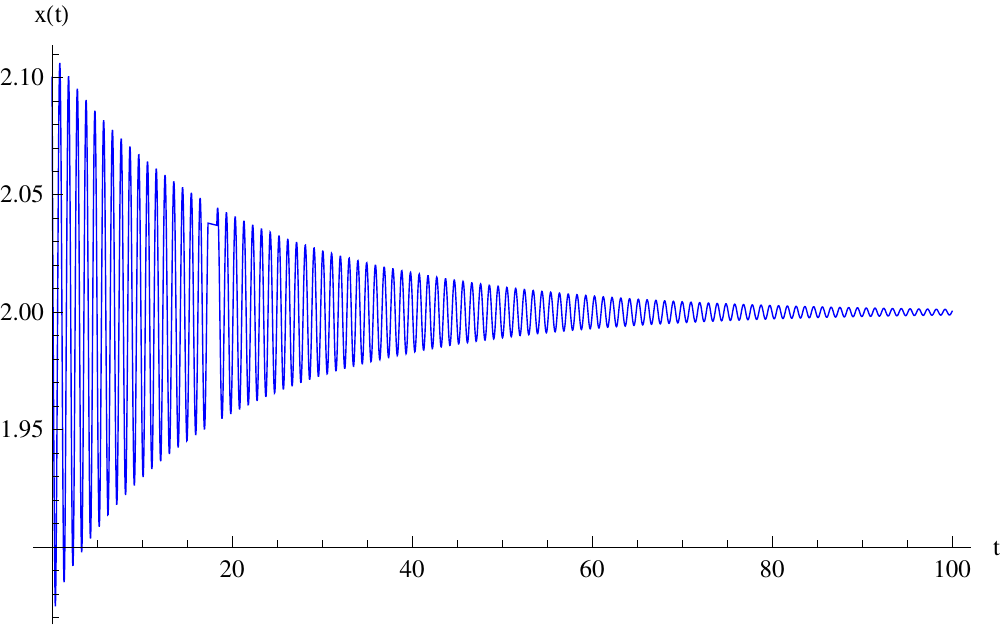}
\caption{The stable solution of (\ref{ex}) with $\alpha=1$ (integer order)}\label{fig1.10}
\end{figure}

\begin{figure}
\centering
\includegraphics[width=0.6\textwidth]{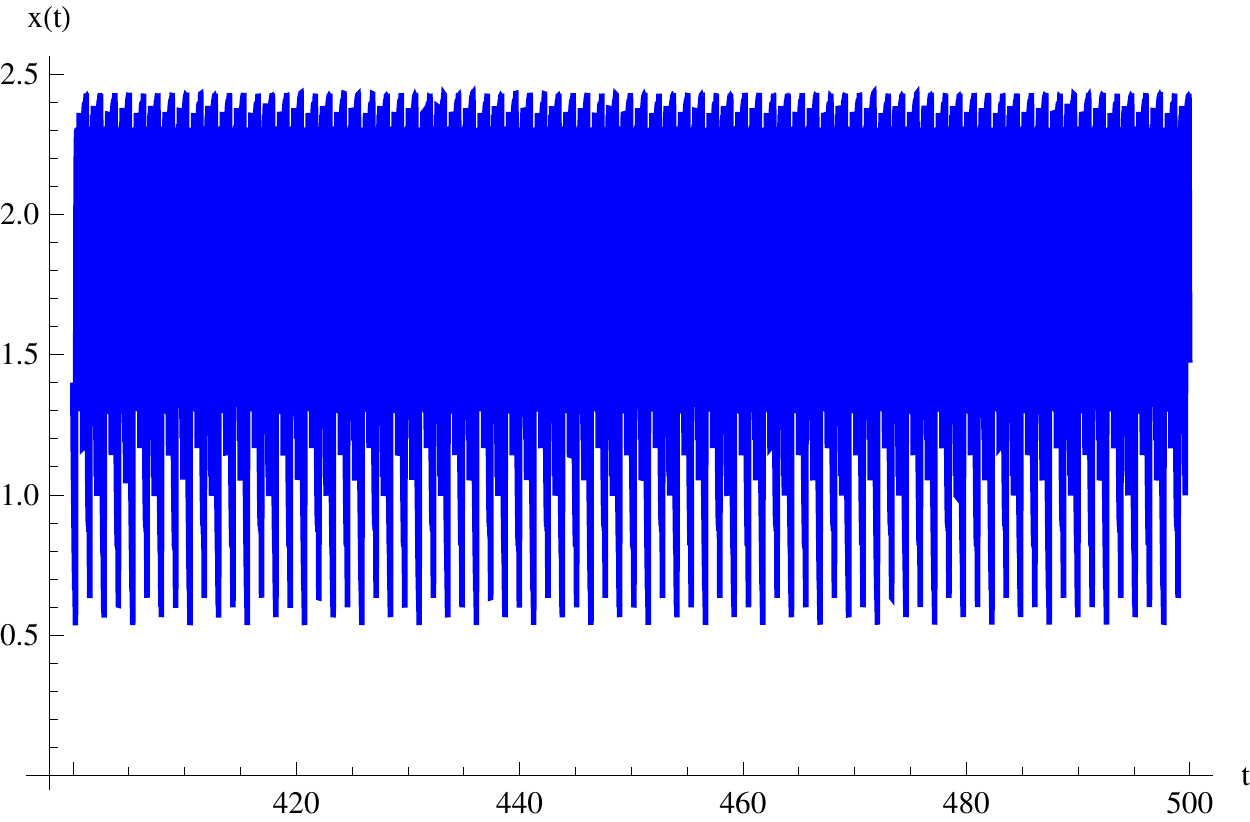}
\caption{Time series of (\ref{ex}) for $\alpha=0.27$}\label{fig1.3}
\end{figure}

\begin{figure}
\centering
\includegraphics[width=0.6\textwidth]{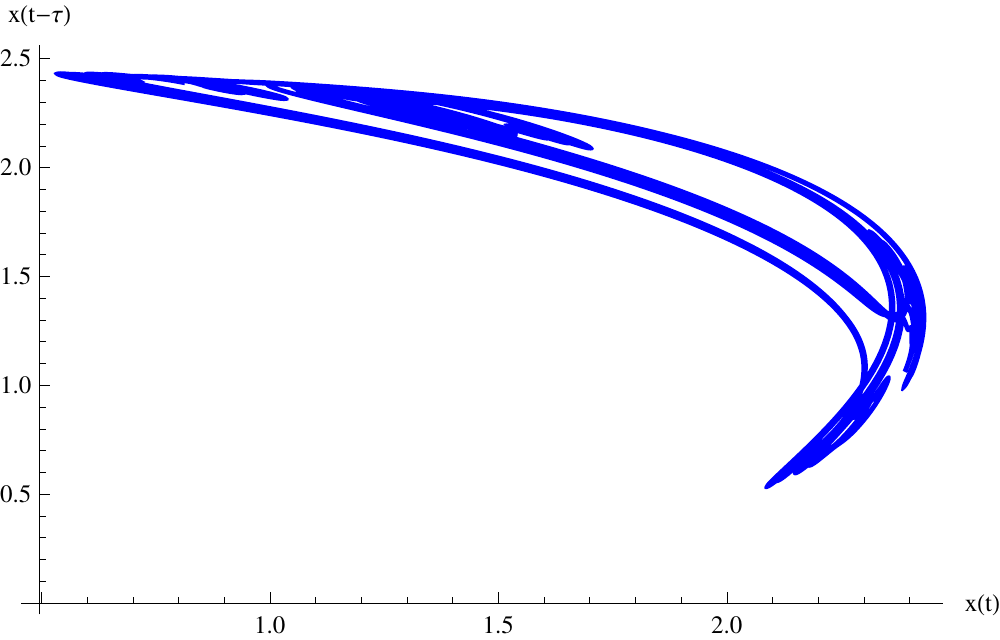}
\caption{The chaotic attractor of (\ref{ex}) at $\alpha=0.27$}\label{fig1.4}
\end{figure}

In this section, we provide the answer to the question in the title of this paper. For the FDEs without delay, the stability of the system at higher value of fractional order $\alpha$ does not get disturbed for the smaller values of the order. As we have discussed in section \ref{sec2}, in the region DS1 of Figure \ref{fig1.1}, the fractional order system may be stable for higher orders and unstable for the lower order values. In the following example, the parameters $a$ and $b$ of the linearized equation near the equilibrium $x^*=2$ belongs to the region DS1. We show that, the system is stable for $\alpha=1$
and show chaos for $\alpha=0.27$. To the best of authors' knowledge, this is first time reported in the literature. It can be observed that the chaotic attractor in this example is very similar to that in the discrete dynamical system viz. Henon map \cite{gallas1993structure}.\\
\textbf{Example} Consider the nonlinear FDDE
\begin{equation}
D^\alpha x(t)=x(t)-x^2(t)+5x(t-\tau)-x(t-\tau)^3,\quad 0<\alpha\leq1. \label{ex}
\end{equation}
It is clear that $x^*=2$ is an equilibrium point.\\
The linearized equation at $x^*$ is \[D^\alpha x(t)=-3x(t)-7x(t-\tau)\].\\
Therefore, $a=-3$ and $b=-7$.\\
Here, $a_0=-4.10117.$\\
Therefore, $b<a_0$ and the stability region for $x^*$ is DS1 (as shown in Figure \ref{fig1.1}).\\
If we fix $\tau=0.31$ then the corresponding critical value of fractional order is $\alpha_*=0.93777.$\\
Therefore, For $\alpha=1$, we get the stable solution as shown in Figure \ref{fig1.10}.\\
Further, for $\alpha=0.27<\alpha_*$, we get a chaotic solution as shown in Figure \ref{fig1.3}. The chaotic attractor in this case is as shown in Figure \ref{fig1.4}. The maximum Lyapunov exponent \cite{kodba2004detecting} in this case is $9.885276.$
\section{Conclusions}\label{sec6}
In this paper, we extended the analysis in \cite{bhalekar2016stability} and revealed the various bifurcation regions in the delay dependent stable region.\\
We showed that, there are three subregions viz. DS1, DS2 and DS3. The region DS3 has analogous bifurcation properties as the non-delayed FDEs whereas the region DS1 shows the exactly contrast properties. The region DS2 has richer dynamics and show various bifurcation scenarios including the stable window. The remarkable acheivement is the example which shows the stable solutions for integer order $\alpha=1$ and chaos for smaller value $\alpha=0.27$.\\
We hope that this manuscript will open new avenues in the research in fractional order delay differential equations.
\section{Acknowledgments}
S. Bhalekar acknowledges the University of Hyderabad for Institute of Eminence-Professional Development Fund (IoE-PDF) by MHRD (F11/9/2019-U3(A)).
D. Gupta thanks University Grants Commission for financial support (No.F.82-44/2020(SA-III)).
\bibliographystyle{elsarticle-num}
\bibliography{paper2.1.bib}
\end{document}